\documentclass[11pt] {article} 
\usepackage[backref]{hyperref}
  \usepackage{amsmath}
    \usepackage{amssymb}
  \usepackage[dvipsnames]{xcolor}   
     \usepackage{pdfsync}

\newtheorem{example}{Example}[section]
\newtheorem{theorem}{Theorem}[section]
\newtheorem{lemma}{Lemma}[section]

\newtheorem{remark}{Remark}[section]

\newcommand{\eqnsection}{
   \renewcommand{\theequation}{\thesection.\arabic{equation}}
   \makeatletter
   \csname @addtoreset\endcsname{equation}{section} 
   \makeatother}

\def \ov{\overline}

\def \be{\begin{equation}}
\def \ee{\end{equation}}
\def \bt{\begin{theorem}} 
\def \et{\end{theorem}}
\def \bl{\begin{lemma}} 
\def \el{\end{lemma}}
\def \bea{\begin{eqnarray}}
\def \eea{\end{eqnarray}}
\def \bas{\begin{eqnarray*}}
\def \eas{\end{eqnarray*}}

\def \al{\alpha}
\def \bb{\beta}
\def \ga{\gamma}
\def\Ga{\Gamma}
\def \de{\delta}
\def \De{\Delta}
\def \ep{\epsilon}

\def \la{\lambda}

\def \Om{\Omega}

\def \vf{\varphi}
\def \si{\sigma}

\def \th{\theta}

\def \ff{\infty}
\def \wh{\widehat}
\def \wt{\widetilde}

\def \FF{{\cal F}}

\def \KK{{\cal K}}

\def \MM{{\cal M}}

\def\b1{\mathbf 1}
\def \({\left(}
\def \){\right)}

\def \da{\downarrow }

\def \nn{\nonumber}
 
\def \Proof{\noindent{\bf Proof $\,$ }}

\def \bc{\begin{center} }
\def \ec{\end{center} }
\def \bs{\begin{slide} }
\def \es{\end{slide} }

\def\square{{\vcenter{\vbox{\hrule height.3pt
        \hbox{\vrule width.3pt height5pt \kern5pt
           \vrule width.3pt}
        \hrule height.3pt}}}}
\def\qed{{\hfill $\square$ \bigskip}}

\eqnsection

\begin{document}

\title{ Local moduli of continuity  for    permanental processes that are zero at  zero.   }

 \author{  Michael B. Marcus\,\, \,\, Jay Rosen \thanks{Research of     Jay Rosen was partially supported by  grants from the Simons Foundation.   }}
\maketitle
 
 \begin{abstract} 
 Let $u(s,t)$ be a continuous potential density of  a symmetric L\'evy process or  diffusion with state space $T$ killed   at $T_{0}$, the first hitting time of $0$,  or   at   $\lambda \wedge T_{0}$, where $\lambda$ is an independent  exponential time. Let
 \[
    f(t)=\int_{T}  u(t,v)\,d\mu(v),
\] 
  where $\mu$ is a finite positive measure on $T$. Let   $X_{\alpha}=\{X_{\alpha}(t),t\in T \}$ be an $\al-$permanental process with kernel  
  \[
  v(s,t)=u(s,t)+f(t).
\] 
 Then when $\lim_{t\to 0}u(t,t)=0$,  \[
  \limsup_{t\downarrow 0}\frac{X_{\alpha}(t )}{u(t,t)\log \log 1/t }\ge 1 ,\qquad  \text{a.s.} \]
and
 \[
  \limsup_{t\downarrow  0}\frac{X_{\alpha}(t )}{u(t,t)\log \log 1/t }\le 1+C_{u,h} ,\qquad  \text{a.s.} \]
where $C_{u,\mu}\le |\mu|$ is a constant that depends on both $u$ and $\mu$, which is given explicitly, and is different in the different examples.    
        
\end{abstract}

\maketitle

   \section{Introduction}\label{sec-1}

Let   $K=\{K(s,t),s,t\in T \}$ be a kernel with the property that for all ${\mathbf t_n} =(t_{1},\ldots,t_{n})$ in $T^n $, the matrix  $\mathbf{ \KK(      t_n)}=\{K(t_{i},t_{j}),i,j\in [1,n] \}$   determines an $n-$dimensional random variable random variable $ (X_{\al}(t_{1}),\ldots, X_{\al}( t_{n}))$ with  Laplace transform,   
\begin{equation}
  E\(e^{-\sum_{i=1}^{n}s_{i}  X_{\al }(t_i) }\) 
= \frac{1}{ |I+ \mathbf{ \KK(      t_n)}S(\mathbf{ s_n })|^{ \al}},   \label{int.1}
\end{equation}
  where $S( \mathbf{        s_n} )$ is a  diagonal matrix with positive entries $\mathbf{ s_n }=( s_{1},\ldots,s_{n} )$, and $\al>0$.     We refer to  $ \mathbf{ \KK(      t_n)}$ as the kernel of $(X_{\al }(t_1),\ldots,X_{\al }(t_n))$.     
It follows from the Kolmogorov Extension Theorem that $\{K(s,t),s,t\in T \}$ determines  a stochastic process which we denote by   $X_\al =\{X_{\al}(t),t\in T \}$ and refer to as an $\al-$permanental process   with kernel $K$.

There are well known examples of $\al-$permanental processes when $\al=k/2$.
Let  $\eta=\{\eta(t);t\in  T\}$  be a mean zero   Gaussian process with covariance $U=\{U(s,t),s,t\in  T\}.$   Let $\{\eta_{i};i=1,\ldots, k\}$ be independent copies of $\eta$. Define,
\begin{equation} \label{1.9mm}
 Y_{k/2}(t)=\sum_{i=1}^{k}\frac{\eta^2_{i}(t)}{2},\qquad t\in T.\end{equation}
The positive stochastic process  $\{Y_{k/2}(t),t\in T \}$ is referred to as a chi--square process of order  $k$  with covariance $U$.
It is well known that for $\{t_1 ,\ldots,t_n\}\subset  T$  
\begin{equation}
  E\(e^{ -\sum_{j=1}^{n}s_{j} Y_{k/2}(t_{j})   }\)
= \frac{1}{ |I+  \ov U(   \mathbf{    t_n})S(\mathbf{        s_n })|^{k/2}}, \label{int.1pp} 
\end{equation}
 where $   \ov U(\mathbf{t_n})$ is the $n\times n$ matrix, $\{
U(t_{j},t_{k}); 1\leq j,k\leq n\}$.  
 Therefore $Y_{k/2}=\{Y_{k/2}(t);t\in T\}$ is a  $k/2$--permanental process with  kernel   $U$.  

 It is important to note that whereas (\ref{int.1pp}) holds for all integers $k$ whenever $U$ is the covariance of a Gaussian process, it may or may not hold when $k/2$
   is replaced by some $\al>0$ that is not a multiple of 1/2. (The question of when 
(\ref{int.1pp}) holds for all $\al>0$ was raised by Paul L\'evy when he asked when does a Gaussian process have infinitely divisible squares. This was answered by 
Griffiths and Bapat. An extensive treatment of these results is given in \cite[Chapter 13]{book}.) Nevertheless,
it follows from Eisenbaum and Kaspi,  \cite[Theorem 3.1]{EK} that  if   $  U $ is 
  the    potential density of a transient Markov process on   $T$ which is   finite for all $s,t\in T$,   then there exist     $\al-$permanental processes  $\{Y_{\al}(t),t\in T \}$ with kernel      $  U$ for all $\al>0$.  I.e., (\ref{int.1pp}) holds with $k/2$ replaced by $\al$  for  all $\al>0.$

  It is at this point that this   becomes interesting because in the   Eisenbaum--Kaspi    theorem   the  kernel    $ U$ is not necessarily  symmetric,  or equivalent to a symmetric  matrix;  (see \cite{MRnonsym}).    When it is not,   the corresponding permanental processes are not   squares of a Gaussian process. They are a new kind of positive stochastic process defined by kernels that need not be symmetric. It is intriguing  to study  sample path properties of these processes.  
  
\medskip  In  \cite{MRLIL} we consider    symmetric    L\'evy processes and   diffusions $Y$,   that are killed at the end of an independent exponential time or the first time they   hit  0, with 
   potential    densities    $u(x,y)$ with respect to some $\sigma$-finite measure on  $T\subseteq R^{1}$.
Let $f$  be  a finite excessive function     for  $ Y$.  Set
 \be u_{  f}(x,y)= u(x,y)+ f(y),\qquad x,y\in  T. 
 \ee

 Under  general smoothness  conditions on     $f$,
  $u(x,y)$ and points  $d\in T$, laws of the iterated logarithm are found for  $X_{k/2} =\{X_{k/2}(t), t\in T \}$, a $k/2-$permanental process with kernel   $ \{u_{  f}(x,y),x,y\in T \}$, of the following form:
For   all integers $k\geq 1$,   
 \be\label{1.4mm}
 \limsup_{x \to 0}\frac{|  X_{k/2}( d+x)-  X_{k/2}  (d)|}{ \left(   2 \sigma^{2}\left(x\right)\log\log 1/x\right)^{1/2}}=    \left( 2 X  _{k/2} (d)\right)^{1/2}, \qquad a.s. , 
 \ee
    where,  
 \be 
   \sigma^2(x)=u(d+x,d+x)+u(x,x)-2u(d+x,x). 
\ee

\medskip In this paper we consider the question of finding the  local modulus of continuity for    $X_{\al}$ at 0 when $X_{\al}  (0)=0$. This, it seems,  is difficult, and, in general,   we can only get different upper   and lower bounds.   On the other hand our results hold  for $X_{\al}$ for all $\al>0$.   

\medskip We consider the following situation:  Let $T$ a be locally compact set with a countable base. 
 Let      $ Y = 
( \Om,   \FF_{t},   Y_t,  \th_{t},   P^x
)$ be a transient Borel right process with state space $T$,  and  continuous strictly positive  potential   densities $u(x,y)$ with respect to some $\si$--finite measure $m$ on $T$.    Let $\mu$ be a finite  positive measure on $T$   with mass $|\mu|=\mu(T)$. 
We call the   function   
\begin{equation}
f(x)= \int_{T} u(y,x)\,d\mu(y)\label{rp.1}
\end{equation}
a left potential for $  Y$. 
 
\medskip The following is \cite[Theorem 6.1]{MRejp}.

\bt\label{theo-borelpr}
Let   $ Y$ be a transient Borel right process with state space $T$,   as described above. Then for any left potential $  f$  for $   Y$, there exists a transient Borel right process   
$  \wt Y \!=\!
(  \Om,   \FF_{t},   \wt Y_t,  \th_{t},   P^x)$ with state space $ S=T\cup \{\ast\}$, where $\ast$ is  an isolated point, such that  $\wt Y$ has potential densities, 
\begin{eqnarray}
&&  \wt    u(x,y)= u(x,y)+ f(y), \hspace{.2 in}x,y\in T
\label{rp.2}\\
&&   \wt    u(\ast,y)=   f(y), \hspace{.2 in}\mbox{ and }\hspace{.2 in} \wt    u(x, \ast)=\wt    u(\ast, \ast) =1, 
\nonumber
\end{eqnarray}
with respect to the measure $\wt m$ on $S$ which is equal to $m$ on $T$ and assigns 
a unit mass to $\ast$.  \et

 Let $\{Z_{t},  t\in S \}$ be the $\al$ permanental process with kernel $\wt    u(x,y)$. 
 We have the following   general lower bound:

\begin{theorem} \label{theo-7.2nn}
 Let      $ Y $ be a strongly symmetric transient Borel right process with state space $(0,1)\subseteq T\subseteq R^{1}$  and  continuous strictly positive  potential   densities $u(x,y)$ with respect to some $\si$--finite measure $m$ on $T$,   as in Theorem \ref{theo-borelpr}.   Assume     that
    $ \lim_{t\downarrow0} u(t,t)=0.$  Then 
   \begin{equation} \label{7.3mmww}
 \limsup_{t\to  0}\frac{Z_t}{u(t,t)\log\log 1/t}\ge 1,\qquad a.s.
\end{equation}
\end{theorem}

\begin{remark}\label{rem-1}{\rm In    Examples 1.1--1.4, $T=R^{1}-\{0\}$ or $(0,\ff)$ and 
\begin{equation}
\wt    u(x,y)= u(x,y)+ f(y), \hspace{.2 in}x,y\in (0,1]\label{si.1}
\end{equation}
has a continuous extension to $x,y\in [0,1]$. We continue to denote this extension by $\wt    u(x,y)$. It follows from the proof of \cite[Lemma 6.2]{MRejp} that this extension is the kernel of an $\al-$permanental process.
   }\end{remark}

We use  \cite[Theorem 1.2]{MRejp}
 to get upper bounds. We repeat it here for the convenience of the reader.

   \begin{theorem} \label{theo-1.1} Let $ X_\al=\{X_\al(t),t\in [0,1] \}$ be an $\al-$permanental process with kernel $ \{\wt u(s,t),s,t\in [0,1] \}$  with  $\wt u(0,0)=0$.      Define the function,     \begin{equation} \label{si.1a}
 \wt \si(s,t)=\(\wt u(s,s)+\wt u(t,t)-2(\wt u(s,t)\wt u(t,s))^{1/2}\)^{1/2}.
\end{equation}
 Let $\vf$  be a   bounded  increasing function  which is regularly varying at zero with positive index such that for all $s,t\in [0,1]$,
 \begin{equation} \label{1.12mm}
 \wt \si(s,t)\le \vf(|s-t|),
\end{equation}  
If, in addition, $\vf^{2}(h)=O(\wt u(h,h))$ at $0$, then
\begin{equation} \label{6.2}
  \limsup_{t\da 0}\frac{X_\al(t )}{\wt u(t,t)\log \log 1/t }\le 1 ,\qquad \text{a.s.} \end{equation}
\end{theorem}

  Note that when
 \begin{equation} \label{1.15mm}
  \wt u(s,t)=u(s,t)+f(t),
\end{equation}
  the denominator in (\ref{6.2}) is larger than the denominator in (\ref{7.3mmww}).

\medskip   Let us now consider Theorem \ref{theo-1.1} with  $\wt u(s,t)$ as in (\ref{si.1})  and 
 \begin{equation} \label{1.17mm}
  f(t):=f_{u,\mu}(t)=\int_{T} u(t,v)\,d\mu(v),\end{equation}
   for some positive finite measure $\mu$.  
By  \cite[Lemma 3.4.3]{book},    $u(t,v)\le u(t,t)\wedge u(v,v)$ so that
 \begin{equation} \label{1.17mms}
 f(t) \le u(t,t)|\mu| 
\end{equation}
In this case we can write (\ref{6.2}) as,
 \begin{equation} \label{6.2ww}
  \limsup_{t\da 0}\frac{X_\al(t )}{  u(t,t)\log \log 1/t }\le 1 +|\mu|,\qquad \text{a.s.} \end{equation}
 
    However, in    the examples we consider, (\ref{6.2ww}) is generally too large. In    the examples we consider we show that,  \begin{equation} \label{1.13} 
  f(t)\sim C_{u,\mu} u(t,t) ,\qquad \text{as $t\da 0$},
\end{equation}
where   $0<C_{u,\mu}\le |\mu|$ is a constant that depends on both $u$ and $  \mu$. In this case   (\ref{6.2}) can be written as,
\begin{equation} \label{6.2q}
  \limsup_{t\da 0}\frac{X_\al(t )}{u(t,t)\log \log 1/t }\le  1+C_{u,h} ,\qquad \text{a.s.} \end{equation}
  
  We consider four examples in which $u(s,t)$  are   continuous potential densities of strongly symmetric Borel right processes killed the first time they hit 0. In the first two examples   the strongly symmetric Borel right  processes are   L\'evy processes,  or exponentially killed L\'evy processes, killed when they hit zero.

\begin{example}\label{ex-1.1}{\rm \label{1.15q}
  Let $Z=\{Z_{t};t\in R^1 \}$ be a real valued symmetric 
L\'evy process with characteristic exponent $\psi$, i.e., 
\begin{equation} \label{1.27nn}
  E\(e^{i\la Z_{t}}\)=e^{-t\psi(\la)},\qquad \forall\,t\ge 0,
\end{equation}
where  $\psi(\la)$ is an even function and,
\begin{equation} \label{1.20nn}
  \frac{1}{ 1+\psi(\la)}\in L^{1}(R^1).
\end{equation} 
    In this example we consider processes $Z$ for which $0$ is recurrent. 
    
     Let $Z' =\{Z'(t);t\in R^1 \}$ be a strongly symmetric  Borel right  process with state space   $T=R^1-\{0\}$ obtained by killing  $Z$ the first time it hits $0$. It follows from \cite[Theorem 4.2.4]{book}  that the process   $Z'$ has  continuous  potential densities with respect to Lebesgue measure restricted to $T$, which we denote by $\Phi  =\{\Phi  (x,y),x,y\in T\}$,  where  
\begin{equation} \label{u01}
\Phi  (x,y)=\frac{(\si^0)^2(x)}{2}+\frac{(\si^0)^2(y)}{2}-\frac{(\si^0)^2(x-y)}{2} ,\qquad x,y\in R^1-\{0\},  
\end{equation}
and 
\be 
   (\si^0)^2(x)=\frac{2}{\pi}\int_{0}^{\ff}\frac{1- \cos\la x}{  \psi(\la)}\,d\la .
\ee 
Note that,  
\begin{equation}
\Phi  (x,x)= (\si^0) ^2(x).\label{rec1.20}
\end{equation}

As we point out in Remark \ref{rem-1} there exists an $\al-$permanental process $\{X_{\al } (x),x\in [0,1]\}$ with kernel,
\be  \label{1.27}
  v(x,y)=\Phi(x,y)+ f_{\Phi,\mu }(y), 
\ee 
where
 \be  \label{3.21a}
   f_{\Phi,\mu }(y) =\int_T \Phi(x,y)  \,d\mu(x)  
\ee 
and   $ |\mu|<\ff$.  

\medskip The next theorem is an application of Theorem \ref{theo-1.1} with $\wt u=v$. 

\begin{theorem}\label{theo-1.3q}  Let $\{X_{\al } (x),x\in [0,1] \}$ be an $\al-$permanental process with kernel $v  (x,y)$.  
 When  $ (\si^0) (t)\in C^1(0,\ff)$   is   a regularly varying   function at zero  with positive index,     
 $X_{\al } (t)$ is continuous on $[0,\de]$. Furthermore, when 

  \begin{equation} \label{1.20}
  \lim _{t\da 0}\frac{t}{ (\si^0) ^{2}(t)}=0,
\end{equation}    \begin{equation}      \limsup_{  t \downarrow 0}\frac{   X_{\al } (t)}{  (\si^0) ^{2}(t) \log\log 1/t }\le 1+\frac{|\mu|}{2} ,   \qquad a.s.,  \label{1.21}
   \end{equation}
 and when  
  \begin{equation} \label{1.22}
  \lim _{t\da 0}\frac{t}{(\si^0)  ^{2}(t)}=\frac{1}{C},\qquad C>0,  
\end{equation}  
    \begin{equation} \label{1.23}    \limsup_{  t \downarrow 0}\frac{   X_{\al } (t)}{  C  t \log\log 1/t }\le     \(1+\frac{|\mu|}{2}\)  +\frac{1}{ 2C}\int _{-\ff}^\ff   ( (\si^0)^2(v))'   \,d \mu(v) ,   \quad a.s. 
   \end{equation}
   \end{theorem}
We give conditions when  $ (\si^0) (t)\in C^1(0,\ff)$ in Remark \ref{rem-c1}.
}\end{example} 

\begin{remark} {\rm 
We show in \cite[Lemma 4.1]{MRLIL} that either,
\begin{equation} \label{1.33e}
  \lim _{t\da 0}\frac{t}{(\si^0)  ^{2}(t)}=0,\qquad\text{or}\qquad\lim _{t\da   0}\frac{t}{(\si^0)  ^{2}(t)}=C,  
\end{equation}  
 for some constant $C>0.$ 
 
\medskip Note also that when $\mu$ is a symmetric measure, since $((\si^0)^2(v))'$ is an odd function, the integral in (\ref{1.23}) is equal to 0. If $\mu$ is concentrated on thee negative half--line it is negative.
   
}\end{remark}

%

 \begin{example} {\rm \label{ex-1.2}  Let $Z_{\bb }=\{Z_{\bb }\(t\);t\in R^1 \}$ be   a symmetric Markov process obtained by killing $Z$ at the end of an independent exponential time with mean $1/\bb$, $\bb>0$. The $\bb-$potential density of $Z$ with respect to Lebesgue measure,   which is the same as the $0-$potential density of $Z_{\bb }$ with respect to Lebesgue measure,  is    
\begin{equation} \label{1.21nn}
 u^{\bb}(x,y)= u^{\bb}(x-y)=\frac{1}{2\pi}\int_{-\ff}^{\ff}\frac{ \cos\la (x-y)}{\bb+ \psi(\la)}\,d\la,\qquad x,y\in R^{1},  
\end{equation}
where $\psi$ is an even function;   see e.g.,   \cite[(4.84)]{book}.  Let 
  \bea 
  (\si^\bb)^2(x-y))&=&u^\bb(x,x)+u^\bb(y,y)-2u^\bb(x,y)\\
  &=&\frac{1}{\pi}\int_{-\ff}^{\ff}\frac{1- \cos\la (x-y)}{\bb+ \psi(\la)}\,d\la,\qquad x,y\in R^{1}.
\eea
  Let 
 \be  \label{3.21bww}
   f _{u^\bb,\mu}(y) =\int_{-\ff}^\ff u^\bb(x,y)   \,d\mu(x), 
\ee 
where $\mu$ is a finite positive    measure with $\mu(\{0 \})=0$.

 Let $Z'_\bb =\{Z'_{\bb}(t);t\in R^1 \}$ be a strongly symmetric  Borel right  process with state space   $T=R^1-\{0\}$ obtained by killing  $Z_\bb$ the first time it hits $0$. The process   $Z_\bb'$ has  continuous  potential densities with respect to Lebesgue measure restricted to $T$, which we denote by $v^{\bb}=\{v^{\bb}(x,y),x,y\in T\}$,  where   
   \bea\label{1.29}
 v^{\bb} (x,y)&=&u^\bb(x,y)-\frac{u^\bb(x,0)u^\bb(0,y)}{u^\bb(0,0)} \label{80.4a}\\
 &=&u^\bb(x-y)-\frac{u^\bb(x)u^\bb(y)}{u^\bb(0)}; \nn
 \eea
 see \cite[(4.165)]{book}. 

\medskip It follows from  Remark \ref{rem-1} that there exists an $\al-$permanental process $\{X_{\al,\bb } (x),x\in [0,1] \}$ with kernel
\be  \label{1.27b}
 v (x,y)= v^\bb(x,y)+  f_{v^\bb,\,\mu}(y),
\ee 
 where,\be  \label{3.21b}
   f _{v^\bb,\mu}(y) =\int_T v^\bb(x,y)   \,d\mu(x) 
\ee 
for   $ \mu  $  as given above.    Since $\mu(\{0 \})=0$ we can write (\ref{3.21b}) as 
 \be  \label{3.21bs}
   f _{v^\bb,\mu}(y) =\int_{-\ff}^\ff v^\bb(x,y)   \,d\mu(x) 
\ee

The next theorem is an application of Theorem \ref{theo-1.1} with $\wt u=v$.
The reader should note that in (\ref{1.21qq}) and (\ref{1.23qq})  we use   $  f _{u^\bb,\mu}$.

 \begin{theorem}\label{theo-1.4}  Let $\{X_{\al } (x),x\in [0,1] \}$ be an $\al-$permanental process with kernel $v  (x,y)$.  When $ (\si^\bb)^2 (t)\in C^1(T)$ and  $|(\si^\bb)^2 (t)'|$ is bounded on $(-\ff,-\De)\cup (\De,\ff)$ for all $ \De>0$ and   $(\si^\bb)^2 (t) $   is     regularly varying   at zero with positive index,  
  $X_{\al,\bb } (t)$ is continuous on $[0,\de]$, for some $\de>0$. Furthermore, when 
  
  \begin{equation} \label{1.20qq}
  \lim _{t\da 0}\frac{t}{ (\si^\bb) ^{2}(t)}=0,
\end{equation}     \begin{equation}      \limsup_{  t \downarrow 0}\frac{   X_{\al } (t)}{  (\si^\bb) ^{2}(t) \log\log 1/t }\le 1+  \frac{ f_{u^\bb,\mu}(0)}{2u^\bb(0)} ,   \qquad a.s.,  \label{1.21qq}
   \end{equation}
and when  
  \begin{equation} \label{1.22qq}
  \lim _{t\da 0}\frac{t}{(\si^\bb)  ^{2}(t)}=\frac{1}{C},\qquad C>0,  
\end{equation}   
  \begin{equation} \label{1.23qq}    \limsup_{  t \downarrow 0}\frac{   X_{\al } (t)}{  C  t \log\log 1/t }\le      1+ \frac{ f_{u^\bb,\mu}(0)}{2u^\bb(0)}   +\frac{1}{2C}\int_{-\ff}^\ff   ( (\si^\bb)^2(v))'   \,d\mu(v) ,   \quad a.s.     \end{equation}
    \end{theorem}

  In Remark    \ref{rem-c1} we give examples of L\'evy processes that satisfy the hypotheses of this theorem. Also  we show in  \cite[Lemma 8.1]{MRLIL}
that for all $\bb\ge 0$,
 \begin{equation} 
  \lim_{t\da 0}\frac{(\si^\bb)  ^{2}(t)}{(\si^{0})^2(t)}
=1,\qquad {a.s}\end{equation} 
 so we can replace   $(\si^\bb)  ^{2}$ by $ (\si^{0})^2$ in   (\ref{1.21qq}). 
 }\end{example}

%
%
\begin{example}{\rm \label{ex-1.3}

    Let  $\mathcal{Z}$ be a  diffusion in $R^1$ that is regular and without traps  and  is symmetric with respect to a $\si$--finite measure  $m$, called the speed measure, which is absolutely continuous with respect to Lebesgue measure.  
(A diffusion is   regular and without traps when
$P^{x}\( T_{y}<\ff\)>0,   \forall x,y\in  R^1$.) We   consider diffusions $\mathcal{Z}$ with  generators of the form
\begin{equation}
\frac{1}{2}a^{2}(x)\frac{d^{2}}{dx^{2}}+c(x)\frac{d}{dx}\label{gendef.1},
\end{equation}
where $a^{2}(x)$ and $c(x)\in C\(R^1\)$ and $a(x)\ne 0$ for any $x\in R^1$.    For details see   \cite[Section 7.3]{RY} and \cite[Chapter 16]{Breiman}.

 Let  $s\in C^{2}(R^1)$ be a scale function for $\cal Z$. The function $s$ is   strictly  increasing and is unique up to linear transformation so that we take $s(0)=0$.  Assume that $\lim_{x\to \ff}s(x)=\ff$. 

 Let  $\ov Z$ be the process with state space  $T=(0,\ff)$ that is obtained by starting  $\cal Z$ in $T$ and then  killing it  the first time it hits 0.  The potential density of $\ov Z$ with respect to the speed measure  $m$
is,
  \begin{equation} \label{diff.5}
  u_{T_{0}}(x,y)=s(x)\wedge s(y),\qquad  x,y>0,
\end{equation}
where   $s\in C^{2}(0,\ff)$  is   strictly positive and  strictly increasing.

 By Remark \ref{rem-1} there exists an $\al-$ permanental process $\{X_{\al } (x),x\in [0,1] \}$ with kernel
\be v (x,y)=  u_{T_{0}}(x,y)+  f_{u_{T_{0}} ,\mu}(y),
\ee
where \be  \label{3.21c}
   f _{ u_{T_{0}},\mu}(y) =\int_T  u_{T_{0}}(x,y)  \,d\mu(x) 
\ee 
and   $ |\mu|<\ff$.  

\medskip The next theorem is an application of Theorem \ref{theo-1.1} where $\wt u=v$.

 \begin{theorem}\label{theo-1.5} Let $\{X_{\al} (t),t\in [0,1] \}$ be an $\al-$permanental process with kernel $v (x,y)$.
 Then
  $X_{\al } (t)$ is continuous on $[0,\de]$    and, 
   \begin{equation} \label{1.42}
   \limsup_{  t \downarrow 0}\frac{   X_{\al} (t)}{    s(t) \log\log 1/t }\le  1+|\mu| ,   \qquad a.s. 
   \end{equation}
 \end{theorem}
 }\end{example}
(Note that we can replace $s(t)$ by $s'(0)t$ in (\ref{1.42}).)
\begin{example}{\rm  
Let    $Z$ be   $\cal Z$  killed  at the end of an independent exponential time with mean $1/\bb$. The potential  density of $Z$ with respect to the speed measure  $m$ is
 \begin{equation} \label{diff.1}
\wt u^\bb (x,y)= \left\{
 \begin{array} {cc}
 p_{\bb}(x)q_{\bb}(y),& \quad x\leq y  
 \\
 q_{\bb}(x)p_{\bb}(y),& \quad y\leq x   
\end{array}  \right. ,
\end{equation}
where $p_{\bb},q_{\bb}\in C^{2}(R^{1})$  are positive and  $p_{\bb}$ is strictly  increasing and  
$q_{\bb}$ is strictly  decreasing.   (For the dependence of  $p_{\bb}$ and $q_{\bb}$ on $\bb$ see  \cite[(4.14)]{book}.)

Let      $Z'$ be the  process with state space  $T=(0,\ff) $ that is obtained by starting   $  Z$ in $T$ and then  killing it  the first time it hits 0. The potential density of $Z'$ is   
\begin{equation}
\wt v^\bb  (x,y)=\wt u^\bb (x,y)-\frac{\wt u^\bb (x,0)\wt u^\bb (0,y)}{\wt u^\bb(0,0)},\qquad  x,y>0, \label{diff.22}  
\end{equation}
i.e.,
 \begin{equation} \label{diff.1v}
\wt v^\bb (x,y)= \left\{
 \begin{array} {cc}
 p_{\bb}(x)q_{\bb}(y)-\displaystyle\frac{p_{\bb}(0)}{q_{\bb}(0)}q_{\bb}(x)q_{\bb}(y),& \quad 0<x\leq y  
 \\\\
 q_{\bb}(x)p_{\bb}(y)-\displaystyle\frac{p_{\bb}(0)}{q_{\bb}(0)}q_{\bb}(x)q_{\bb}(y),& \quad 0<y\leq x.   
\end{array}  \right. 
\end{equation}

As we point out in Remark \ref{rem-1} there exists an $\al-$permanental process $\{X_{\al } (x),x\in [0,1] \}$ with kernel
\be v(x,y)=  \wt v^\bb (x,y)+ f_{\wt v^\bb,  \mu}(y),\label{xkko}
\ee
 where
 \be  \label{3.21d}
 f_{\wt v^\bb,  \mu}(y)=\int_T \wt v^\bb (x,y) \,d\mu(x),
\ee 
and   $ |\mu|<\ff$.  

\medskip The next theorem is an application of Theorem \ref{theo-1.1} where $\wt  u=v$.

  \begin{theorem}\label{theo-1.6} Let $\{X_{\al} (t),t\in [0,1] \}$ be an $\al-$permanental process with kernel $v  (x,y)$.  
 Then
  $X_{\al } (t)$ is continuous on $[0,\de]$    
    and,  
   \begin{equation} \label{1.49}
   \limsup_{  t \downarrow 0}\frac{   X_{\al } (t)}{   \rho_{\bb}(0) t \log\log 1/t }\le   \(1+\int_{0}^\ff \frac{q_{\bb}(y)}{q_{\bb}(0)}  \,d\mu(y)\),    \qquad a.s. 
   \end{equation}
   where,
   \begin{equation} 
  \rho_{\bb}(0)= {p_{\bb}'(0)}{q_{\bb}(0)}- {q_{\bb}'(0)}{p_{\bb}(0)}.
\end{equation}
 \end{theorem}
 }\end{example}

\section{Proof  of Theorems \ref{theo-7.2nn}   }

  We make  several observations that are   used in the proof of Theorem \ref{theo-7.2nn}.  Let   $U $ be a non-singular $n\times n$ matrix. Let  $U^{-1}$   denote the inverse of $U$ and $U^{j,k}$   denote the elements of $U^{-1}$. 
Let  $U_{f}$ be the $(n+1)\times (n+1)$ matrix   
  \bea \hspace{-.3 in} U_{f}&=&\left (
\begin{array}{ cccc } 1 &f(1 )&\ldots&f(n )  \\
1   &U_{1 ,1 }+f(1 )&\ldots&U_{1 ,n }+f(n )  \\
\vdots& \vdots &\ddots &\vdots  \\
1   &U_{n ,1 }+f(1 )&\ldots&U_{n ,n+ }+f(n )  
\end{array}\right ).\label{19.39}
   \eea 
where $f$ is any function. One can check that 
\be  U_{f}^{-1} =\left (
\begin{array}{ cccc  }  1+\rho &-\sum_{i=1}^{n}f(i)U ^{i ,1 }   &\dots&-\sum_{i=1}^{n}f(i)U ^{i ,n }  \\
- \sum_{j=1}^{n}U ^{1 ,j } & U ^{1 ,1 } & \dots &  U ^{ 1 ,n }    \\
\vdots&\vdots& \ddots&\vdots  \\
- \sum_{j=1}^{n}U ^{n ,j }& U ^{n ,1 }&  \dots & U ^{n ,n }   \end{array}\right ),\label{19.40q}
   \ee 
   where
   \begin{equation}
\rho =\sum_{j=1}^{n}\sum_{i=1}^{n}f(i)U ^{i ,j} .\label{19.10vq}
\end{equation}
We note the following remarkable property of $U_f$ that we use,   \begin{equation} \label{7.7}
 (U_f^{-1})^{i,i}= U  ^{i-1,i-1},\qquad i=2,\ldots,n+1.
\end{equation}


Let $\xi_\al$ be a gamma random variable with probability density function,
\begin{equation} 
 \frac{x^{\al-1}e^{-x}}{\Ga(\al)}.
\end{equation}
The following Theorem is  \cite[Theorem 1.1]{MRun}:
\begin{theorem} \label{theo-7.3} Let $X=(X_0,X_1,\ldots,X_n)$ be an $\al-$permanental vector with non-singular kernel $K$. Let $A=K^{-1}$ and denote the  diagonal entries of $A$ by $(a_0,a_1,\ldots,a_n)$. Then there exists a coupling between $X$ and an $n-$tuple $(\xi^{(0)}_\al,\xi^{(1)}_\al,\ldots,\xi_\al^{(n)})$ of independent identically distributed. copies of $\xi_{\al}$ such that,
\begin{equation} 
 X\ge  \(a_0^{-1}\xi^{(0)}_\al,a_{1 }^{-1}\xi^{(1)}_\al,\ldots,a_{n }^{-1}\xi_\al^{(n)} \),\qquad a.s.
\end{equation}

\end{theorem}

 \noindent\textbf{Proof of Theorem \ref{theo-7.2nn} }  Let   $\{t_j=\th^{j}\}_{j=1}^{n}$ and let   $X=(X (*),X(t_1),\ldots ,X(t_n)$ be an $\al-$permanental vector with non-singular kernel $K$
  \bea \hspace{-.3 in} K&=&\left (
\begin{array}{ cccc } 1 &f(t_{1}  )&\ldots&f(t_n )  \\
1   &u(t_1 ,t_1 )+f(t_1 )&\ldots&u(t_1 ,t_n )+f(t_n )  \\
\vdots& \vdots &\ddots &\vdots  \\
1   & u(t_ n ,t_1 )+f(t_1 )&\ldots&u(t_ n ,t_n)+f(t_n )
\end{array}\right ).\label{19.39n}
   \eea 
where $f$ is given in (\ref{rp.1}). $A=K^{-1}$ is an $M $ matrix. Denote its   diagonal entries by $(a_0,a_1,\ldots,a_n)$. Let $\MM=\{ u(t_i,t_j)\}_{i,j=1}^{n}$.
 It follows from   (\ref{7.7}) that \begin{equation} 
 a_i=\MM^{i,i},\qquad i=1,\ldots,n.
\end{equation}
Let $\ov\MM$ be the matrix with elements,
\begin{equation} \label{mmm.1}
 \ov\MM_{i,j}=\frac{\MM_{i,j}}{\MM^{1/2}_{i,i}\MM^{1/2}_{j,j}},\qquad i,j=1,\ldots,n.
\end{equation}
We have,
\begin{equation} 
 \ov\MM=D\MM D,
\end{equation} 
where $D$ is a diagonal matrix with diagonal elements $\{1/\MM^{1/2}_{i,i}\}_{i=1}^{n}$. 
Therefore,\begin{equation} 
 \MM^{-1}=D \ov \MM^{-1} D .
\end{equation}
Let $  A= \MM^{-1}$ and $\ov A=\ov\MM^{-1}$.   We have, 
\begin{equation} \label{mm.3}
a_{i}=    \MM^{i,i}=A_{i,i} =\(\MM^{-1}\)_{i,i}=\(D \ov \MM^{-1} D \)_{i,i}     = \frac{\ov A_{i,i} }{\MM_{i,i}} ,\qquad i=1,\ldots,n.
\end{equation}
We see by  \cite[Lemma 5.1]{MRun} that,
\begin{equation} 
 \ov A_{i,i}\le \frac1{\ov \MM_{i,i}-\max_{j\ne i}\ov \MM_{i,j}}. 
\end{equation} 
Note that it follows from   the symmetry of $Y$ and \cite[Lemma 3.4.3]{book}
that when $i\ne j$,
\begin{equation} 
 \ov \MM_{i,j}\le \frac{\MM^{1/2}_{i,i}}{ \MM^{1/2}_{j,j}}\wedge\frac{\MM^{1/2}_{j,j}}{\MM^{1/2}_{i,i} }\le \th.
\end{equation}
Therefore,   since $\ov \MM_{i,i}=1$,
\begin{equation} 
 \ov A_{i,i}\le \frac1{1-\th}, 
\end{equation} 
and   by (\ref{mm.3}),  
\begin{equation} \label{noted}
 a_{i}\le \frac{1+2\th}{\MM_{i,i}}=\frac{1+2\th}{u(\th^{ i},\th^{i})},\qquad i=1,\ldots,n.
\end{equation}

It follows from Theorem \ref{theo-7.3} that
\begin{equation}
\(X_{t_{1}}, X_{t_{2}}, \cdots, X_{t_{n}} \)\geq  \(a_{1}^{-1}\xi^{(1)}_\al,a_{2}^{-1}\xi^{(2)}_\al, \ldots, a_{n}^{-1}\xi^{(n)}_\al, \)
\end{equation}
so that
\begin{equation}
\(a_{1} X_{t_{1}}, a_{2} X_{t_{2}}, \cdots, a_{n}X_{t_{n}} \)\geq  \( \xi^{(1)}_\al, \xi^{(2)}_\al, \ldots,  \xi^{(n)}_\al \).\label{mm.5}
\end{equation}
By (\ref{noted}) we see that,  
\begin{equation}
\( \frac{X_{t_{1}}}{u(t_{1}, t_{1})} ,  \frac{X_{t_{2}}}{u(t_{2}, t_{2})}, \cdots, \frac{X_{t_{n}}}{u(t_{n}, t_{n})} \)\geq \frac{1}{1+2\th}  \( \xi^{(1)}_\al, \xi^{(2)}_\al, \ldots,  \xi^{(n)}_\al \).\label{mm.5q}
\end{equation}
 Consequently,  
\be  
   \max_{2\leq i\leq n}    \frac{X_{\th^{i}}}{ u(\th^{i},\th^{i} )\log i}   \ge \frac{1}{1+2\th} \max_{2\leq i\leq n}    \frac{ \xi^{(i)}_\al}{\log i} ,\qquad a.s. 
\ee 
 It now follows from  \cite[Lemma 7.2]{MRejp} that for any $\al$ and $\ep>0$ there exists an $l=l (\ep)$ such that,
\begin{equation} \label{2.21}
  P\(\max_{l\le i \le n}  \frac{X_{\th^{i}}}{ u(\th^{i },\th^{i} )\log i}> \frac{1-\ep}{1+2\th}\)\ge 1-\frac{l+1}{n+1}.
\end{equation} 
  Now let $\{Z_{t},  t\in S \}$ be an     $\al-$permanental process with kernel $\wt    u(x,y)$ in (\ref{rp.2}). It follows from (\ref{2.21}) that,
\begin{equation} \label{2.21q}
  P\(\max_{l\le i \le n}  \frac{ Z_{\th^{i}}}{ u(\th^{i },\th^{i} )\log i}> \frac{1-\ep}{1+2\th}\)\ge 1-\frac{l+1}{n+1}.
\end{equation} 
Since we can take $\th$ arbitrarily small we get  (\ref{7.3mmww}).\qed

\medskip The next lemma breaks down the condition in (\ref{1.12mm}) and is used in everything that follows to show that (\ref{1.12mm}) is satisfied. 

\begin{lemma}\label{lem-key}   Let $\wt u$   be as given in (\ref{rp.2}) and assume that $u$ symmetric. Let  $ \wt \si   $ be as given in  (\ref{si.1a}).  If 
\begin{equation} \label{2.14}
   u ( s,s) +    u( t,t)-    2u ( s,t) \le C u(|t-s|,|t-s|),
\end{equation}
and \begin{equation} \label{2.15}
  \big|f(s)-    f(t)\big|\le C u(|t-s|,|t-s|) 
\end{equation} 
for all $s,t\in [0,\de]$,
then  
\begin{equation} \label{6.18}
  \wt \si  ^2(s,t)\le 2C u(|t-s|,|t-s|), \end{equation}
for all $s,t\in [0,\de]$.
\end{lemma}

\Proof Consider   $\wt \si   ^2(s,t)$. Using the fact that $u(s,t)$ is symmetric we can
  write this as,
 \bea \label{3.107qqs}
  \wt  \si  ^2(s,t)  &=&   \wt u ( s,s) +    \wt u( t,t)-    \wt u ( s,t) -   \wt u( t,s)\\
 &&\qquad    + \( \wt u ^{1/2}( s,t)-       \wt u ^{1/2} ( t,s)\)^2\nn\\&=& u ( s,s) +    u( t,t)-    2u ( s,t)  + \( \wt u^{1/2}(s,t)- \wt u^{1/2}(t,s)\)^2. \nn  \eea
 Once again, using the fact that $u(s,t)$ is symmetric we   have,
 \be  \label{3.107qq2}
\(\wt u^{1/2}(s,t)- \wt u^{1/2}(t,s)\)^2 \le   \big| \wt u (s,t)- \wt u (t,s)\big|\le \big|f(s)-    f(t)\big|. 
\ee 
Now using (\ref{2.14}) and (\ref{2.15})    to bound the last line of (\ref{3.107qqs}) we get (\ref{6.18}).\qed

  We use Lemma \ref{lem-key} to get a more usable version of Theorem \ref{theo-1.1}.

 \begin{theorem} \label{theo-2.2}  
  Let $ X_\al=\{X_\al(t),t\in [0,1] \}$ be an $\al-$permanental process with kernel $ \{\wt u(s,t),s,t\in [0,1] \}$    where $\wt u(0,0)=0$.
  If $u(v,v)$
is   regularly varying   at zero  with positive index and (\ref{2.14}) and (\ref{2.15}) hold   
 then
\begin{equation} \label{6.2ss5}
  \limsup_{t\da 0}\frac{X_\al(t )}{(u(t,t)+f(t))\log \log 1/t }\le 1 ,\qquad \text{a.s.} \end{equation}
\end{theorem}

\Proof If $u(v,v)$
is   regularly varying at zero   with positive index it is asymptotic   to an increasing function at 
zero. Denote this function by $\wh\vf^{2}(v)$.  
 It follows from (\ref{6.18}) that  
\begin{equation} \label{nn}
 \wt  \si  ^2(s,t)\le  C' \wh\vf(|t-s|),\qquad  
\end{equation}
and, since  $\wh \vf^{2}(v)=O(  u(v,v))$ at 0, we have   $\wh \vf^{2}(v)=O( \wt u(v,v))$ at zero.
Consequently,  (\ref{6.2ss5}) follows from (\ref{6.2}).\qed

    We proceed to prove Theorems \ref{theo-1.3q}--\ref{theo-1.6} by establishing   (\ref{2.14}) and (\ref{2.15}).

\section{Proof of Theorem \ref{theo-1.3q}}
  By definition, 
 \begin{equation} \label{6.14q}
   \Phi ( s,s) +     \Phi( t,t)-    2 \Phi ( s,t) =   \Phi(|t-s|,|t-s|)= (\si^0)^2(t-s),
\end{equation}
first for all $ s,t\in R^1-\{0\} $. Since $\Phi(s,t)$ is continuous for all $ s,t\in R^1 $ we   take it to be defined and continuous on $R^1$.   This gives (\ref{2.14}).

\medskip We have the following remarkable regularity property of $(\si^0)^2$ that is proved at the end of this section. 

\begin{lemma} \label{lem-3.1}
   \begin{equation} \label{u01mm}
 |(\si^0)^2(y) - (\si^0)^2(x-y)|\le(\si^0)^2(x) ,\qquad x,y\in R^1,  
\end{equation}

\end{lemma}

Using this we show  that,
\be 
|f_{\Phi,\mu}( s)-  f_{\Phi,\mu}(t)|\le (\si^0)^2(t-s) |\mu|,  \quad s,t\in R^1. \label{zip}\ee
It follows from Lemma \ref{lem-3.1}
that
 \begin{equation} \label{3.22}
  |(\si^0)^2(s-v)-(\si^0)^2(t-v) |\le  (\si^0)^2(t-s)  .
\end{equation}
Therefore,  
\bea 
\lefteqn{   |f_{\Phi,\mu}( s)-  f_{\Phi,\mu}(t)|\le \int _T  |\Phi (s,v)-\Phi(t,v)| \,d\mu(v)  }\\
   & &\le\frac{1}{2} \int_T \(|(\si^0)^2(t)-(\si^0)^2(s) |+|(\si^0)^2(t-v)-(\si^0)^2(s-v) |\)d\mu(v) \nn\\
   &  &  \le\frac{1}{2}\int_T  (\si^0)^2(t-s) d\mu(v) = (\si^0)^2(t-s) |\mu|\nn.
\eea
   The relationships in (\ref{6.14q}) and (\ref{zip}) and Theorem \ref{theo-2.2} show that \begin{equation}      \limsup_{  t \downarrow 0}\frac{   X_{\al } (t)}{  \((\si^0) ^{2}(t)+f_{\Phi,\mu}(t) \)\log\log 1/t }\le 1  ,   \qquad a.s.  \label{1.21s}
   \end{equation}

\medskip We now get an estimate in the form of (\ref{1.13}) so that we can write (\ref{1.21s}) in the form of (\ref{6.2q}). 
Using   Lemma  \ref{lem-3.1} we have, 
 \bea \label{3.9a}
        f_{\Phi,\mu}(t) &=& \int _ T  \Phi(t,v) d\mu(v)= \frac{1}{2}\int_ T \( (\si^0)^2(t)+(\si^0)^2(v)-(\si^0)^2(t-v)\) d\mu(v)\nn\\
   &= & \frac{(\si^0)^2(t)}2{} |\mu|+\frac{1}{2}\int_T\(  (\si^0)^2(v)-(\si^0)^2(t-v)\) d\mu(v).
   \eea
(Here we use the fact that $(\si^0)^2$ is an even function and we extend it  to $T $.) Consequently 
   \begin{equation}
  \frac{f_{\Phi,\mu}(t)}{(\si^0)^2(t)}= \frac{|\mu|}{2} +\frac{1}{2} \int_ T \frac{\(  (\si^0)^2(v)-(\si^0)^2(v-t)\)}{(\si^0)^2(t)} \,d\mu(v). \label{3.223}
   \end{equation}
  By (\ref{3.22})
   \begin{equation}
   \frac{|  (\si^0)^2(v)-(\si^0)^2(v-t)|}{(\si^0)^2(t)}\leq 1.\label{3.22a}
   \end{equation}
 Therefore, 
 using   the Dominated Convergence Theorem we see that, 
 \begin{equation}
  \lim_{t\da 0}\frac{f_{\Phi,\mu}(t)}{(\si^0)^2(t)}= { \frac{|\mu|}{2}} + \frac{1}{2}\int_T\lim_{t\da 0} \frac{\(  (\si^0)^2(v)-(\si^0)^2(v-t)\)}{(\si^0)^2(t)} \,d\mu(v). \label{3.2ujq}
   \end{equation}
Note that,
\bea
\lim_{t\da 0}\frac{   (\si^0)^2(v)-(\si^0)^2(v-t) }{(\si^0)^2(t)} &=&\lim_{t\da 0}  \frac{ (\si^0)^2(v)-(\si^0)^2(v-t)}{t}  \,{\frac{t}{(\si^0)^2(t)}}\nn\\
&=&((\si^0)^2(v))'\lim_{t\da 0}  \frac{t}{(\si^0)^2(t)}.
\eea
 Consequently, the integral in (\ref{3.2ujq}) is equal to 
\begin{equation} 
 \frac{1}2{}  \int_T   ( (\si^0)^2(v))'\lim _{t\da 0}\frac{t}{ (\si^0) ^{2}(t)}   \,d\mu(v).
\end{equation}
Using this and (\ref{3.2ujq}) we get (\ref{1.21}) when (\ref{1.20}) holds and (\ref{1.23}) when (\ref{1.22}) holds.\qed

    \noindent \textbf{Proof of Lemma \ref{lem-3.1}}
   Since  $ \Phi(x,y)$ is a continuous potential density  of a strongly symmetric Borel right process with state space $T= R^1-\{0\}$ it follows from \cite[Lemma 3.4.3]{book} that  
 \begin{equation}
 \Phi(x,y)\leq  \Phi(x,x)\wedge  \Phi(y,y)= (\si^0)^2(x)\wedge(\si^0)^2(y), \quad x,y\in R^1-\{0\}. 
\end{equation}
 Therefore, it follows from this and (\ref{u01}) that,
 \begin{equation} \label{u01qmm}
(\si^0)^2(x)+ (\si^0)^2(y) -(\si^0)^2(x-y)\le2(\si^0)^2(x) ,\quad x,y\in R^1-\{0\},  
\end{equation}
 Since $\si^0$  is continuous on $R^1$ we  have, 
 \begin{equation} \label{u01qmmb}
(\si^0)^2(x)+ (\si^0)^2(y) -(\si^0)^2(x-y)\le2(\si^0)^2(x) ,\quad x,y\in R^1. 
\end{equation}
Or equivalently,\begin{equation} \label{u01q}
 (\si^0)^2(y) -(\si^0)^2(x-y)\le(\si^0)^2(x) ,\qquad x,y\in R^1.  
\end{equation}
Furthermore, since $\Phi(x,y)\geq 0$ it follows from (\ref{u01}) that 
 \begin{equation}
(\si^0)^2(x-y)-(\si^0)^2(y)\le  (\si^0)^2(x) . \label{3.22f}
 \end{equation}
Combining this with (\ref{u01q}) we have (\ref{u01mm}).\qed

 \begin{remark}{\rm  It is interesting to note that since $((\si^\bb)^2 )'(t)$ is an odd function whenever $\mu$ is a symmetric measure the integral in (\ref{1.23}) is equal to 0 and we get 
  \begin{equation} \label{1.2qmm}    \limsup_{  t \downarrow 0}\frac{   X_{\al } (t)}{  C  t \log\log 1/t }\le      1+ \frac{ |\mu|}{2 }    ,   \qquad a.s. 
   \end{equation}

}\end{remark}

 \begin{example} {\rm   Suppose  that $\mu$ is concentrated on the negative half line and  \begin{equation} \label{1.22}
  {(\si^0)  ^{2}(t)}= {C|t|}.  
\end{equation}  
 Then 
 \begin{equation}  \frac{1}{ 2C}\int _{-\ff}^\ff   ( (\si^0)^2(v))'   \,d \mu(v)=- \frac{ |\mu|}{2 }. 
   \end{equation}
 Therefore, it follows from  (\ref{1.23}) and Theorem \ref{theo-7.2nn}  that,
\begin{equation} \label{4.26q}
   \limsup_{  t \downarrow 0}\frac{   X_{\al} (t)}{C  t \log\log 1/t }=1,     \qquad a.s.  
   \end{equation}

}\end{example}

  \section{Proof of Theorem \ref{theo-1.4}}
   We have, 
  \bea
  (\si^\bb)^2(x-y))&=&u^\bb(x,x)+u^\bb(y,y)-2u^\bb(x,y)\\
  &=&2\(u^\bb(0) -u^\bb(x-y)\).\nn
\eea
Define,    \begin{equation} 
  (\si_v^\bb)^2(x,y))=v^\bb(x,x)+v^\bb(y,y)-2v^\bb(x,y).
\end{equation}
Using   (\ref{1.29})  we see that,
\begin{equation} \label{3.xxq}
  (\si_v^\bb)^2(x,y)\le (\si^\bb)^2(x-y) .
\end{equation}
In addition,
\bea \label{4.5}
  v^\bb(x,x)&=&u^\bb(0)-\frac{(u^\bb(x))^2}{u^\bb(0)}=\frac{(u^\bb(0))^2- (u^\bb(x))^2} {u^\bb(0)}\\
   &=&\frac{(\si^\bb)^2(x)}{2u^\bb(0)}\(u^\bb(0)+u^\bb(x)\),\nn
\eea
so that,
\begin{equation} \label{4.6}
  v^\bb(x,x)\sim  { (\si^\bb)^2(x)},\qquad\text{as $x\da 0.$}\end{equation}
  It follows from  (\ref{3.xxq}) and  (\ref{4.5})   that, 
  \begin{equation} \label{4.6xxq}
  (\si_v^\bb)^2(x,y)\leq Cv^\bb(x-y,x-y),\qquad\text{for all  $x,y \in [0,\de].$}\end{equation}
 This gives (\ref{2.14}) for this example.
 
\medskip  To obtain (\ref{2.15}) we use following lemma which is given in  \cite[Lemma 7.4.2]{book}: 

\begin{lemma}  Let $\{G(t),t\in R^{1} \}$   be a mean zero stationary Gaussian process  with   continuous  covariance $u(s,t)=u(t-s)$, that is the potential density of a    strongly symmetric Borel right process.  Set,  
\begin{equation} 
  \si^2(t)=E(G(t)-G(0))^2=2\(u(0)-u(t)\).
\end{equation}
Then,
\begin{equation} \label{4.12a}
  |\si^2(t)-\si^2(s)|\le \frac{2u(0)}{u(t)+u(s)}\si^2(t-s).
  \end{equation}
\end{lemma}

 Since $\{u^\bb(x,y)$, $x,y \in R^1  \}$ is the potential density  of a strongly symmetric  Borel right  process  with state pace $R^{1}$ it follows from \cite[Lemma 3.3.3]{book} that  it is positive definite and therefore it is the covariance of a  stationary Gaussian process.   Applying (\ref{4.12a}) to $u^\bb$ and $\si^\bb$ we have,
\begin{equation} \label{4.12w}
  \Big |(\si^\bb)^2(t)- (\si^\bb)^2(s)\Big |\le \frac{2u^\bb(0)}{u^\bb(t)+u^\bb(s)} (\si^\bb)^2(t-s).
  \end{equation}
The potential of $Z_\bb$ is given in (\ref{3.21bww}). 
 We take $\mu(\{0 \})=0$ so that  potential of $Z'_\bb$ is 
 \be  \label{3.21bwwq}
   f _{v ^\bb,\mu}(y) =\int_{-\ff}^\ff  v^\bb(z,y)   \,d\mu(z). 
\ee 
Consequently,
\be  \label{4.12mm}
   \big|f_{v^\bb,\mu}(x)-f_{v^\bb,\mu}(y ) \big| \le   { \big| f_{u^\bb,\mu}(x)-  f_{u^\bb,\mu}(y ) \big|} + \frac{ f_{u^\bb,\mu}(0)}{2u^\bb(0)} \big| (\si^\bb)^2(x)-(\si^\bb)^2(y)\big|.
   \ee
It follows from (\ref{4.12w}) that,
\begin{equation} \label{4.13q}
 \frac{ f_{u^\bb,\mu}(0)}{2u^\bb(0)} \big| (\si^\bb)^2(x)-(\si^\bb)^2(y)\big|\le  \frac{ f_{u^\bb,\mu}(0)}{   (u^\bb(x)+u^\bb(y))}  (\si^\bb)^2(x-y).
\end{equation} 
  Since $u^\bb(x)$ is continuous, $u^\bb(x)>0$,  for all    $x$ sufficiently small. Therefore, for  $\De$ sufficiently small when $x,y$ in $[0,\De]$ we have,
\begin{equation} \label{4.13qq}
 \frac{ f_{u^\bb,\mu}(0)}{2u^\bb(0)} \big| (\si^\bb)^2(x)-(\si^\bb)^2(y)\big|\le  C_\De  (\si^\bb)^2(x-y),
\end{equation}
for some finite constant $C_\De.$

We now consider the first term after the equal sign  in (\ref{4.12mm}).  We have,
 \begin{eqnarray}
 &&| f_{u^\bb,\mu}( x)-   f_{u^\bb,\mu}(y)|
 \label{slid.11}
 \\
 &&\qquad=\bigg|\int_{-\ff}^\ff u^\bb(z-x) \,d\mu(z)-\int  _{-\ff}^\ff u^\bb( z-y) \,d\mu(z)\bigg|
 \nonumber \\
 &&\qquad=\bigg|\int_{-\ff}^\ff\( \(u^\bb(0)-u^\bb(z-y)\) -\(u^\bb(0)-u^\bb( z-x)\)\)  \,d\mu(z)\bigg|
 \nonumber \\
 &&\qquad \le \frac{1}{2}\int  _{-\ff}^\ff \Big|  (\si^\bb)^2(z-y)-(\si^\bb)^2(z-x)\Big| \,d\mu(z).
 \nonumber
 \end{eqnarray} 
 It follows from (\ref{4.12w})  that  for $\De$ sufficiently small,
  \be  \label{4.22}
  \(\int _0^\De+\int_{-\De}^0 \)\( (\si^\bb)^2(z-y)-(\si^\bb)^2(z-x)\) \,d\mu(z) \le C'_\De(\si^\bb)^2(x-y), \ee
 for some finite constant $C'_\De.$
 
 In addition
 \bea  
  &&\nn H_\De(x-y):=\(\int _{-\ff}^{-  \De}+\int_{\De}^\ff \)\( (\si^\bb)^2(z-y)-(\si^\bb)^2(z-x)\) \,d\mu(z) \\&&\qquad \le \max_{ |z|\ge \De} (\si^\bb)^2(z)'\mu(\{(-\ff,-\De)\cup \mu(\De,\ff)\})|x-y|
   , \eea
and,
\be   
  \frac{H_\De(x-y)}{(\si^\bb)^2(x-y)} \le \max_{ |z|\ge \De} (\si^\bb)^2(z)'\mu(\{(-\ff,-\De)\cup \mu(\De,\ff)\})\frac{|x-y|}{(\si^\bb)^2(x-y)}.
    \ee 

Using this (\ref{4.22}),  (\ref{4.13qq}) and (\ref{1.33e}), we get,
\be  \lim_{x,y\da 0}\frac{ \big|  f_{v^\bb,\mu}(x)-  f_{v^\bb,\mu}(y)\big|}{(\si^\bb)^2(x-y)}\label{4.10mm}  \le   C,\nn   
\ee   
for some constant $C$.
 Together with (\ref{4.6}) this gives us the version of (\ref{2.15}) we need for this theorem.
 
\medskip We  now get an estimate for $ f_{v^\bb,\mu} $ at 0.
  \bea  \label{newq}
  f_{v^\bb,\mu}(x)
   &=&\int_{-\ff}^\ff u^\bb(x,y) \,d\mu(y)-\frac{u^\bb(0,x)}{u^\bb(0,0)}\int _{-\ff}^\ff u^\bb(0,y) \,d\mu(y)
  \\
  &=& \nn  {  f_{u^\bb,\mu}(x)} -\frac{u^\bb(0,x)  f_{u^\bb,\mu}(0)}{ u^\bb(0,0)}\\
  &=&\nn   { f_{u^\bb,\mu}(x)} -  {  f_{u^\bb,\mu}(0)} + \frac{   f_{u^\bb,\mu}(0)}{ u^\bb(0,0)}\( u^\bb(0,0)- u^\bb(0,x) \)\nn\\
   &=& {  f_{u^\bb,\mu}(x)-   f_{u^\bb,\mu}(0)} + \frac{(\si^\bb)^2(x) f_{u^\bb,\mu}(0)}{2u^\bb(0,0)} \nn.
\eea
 In addition,
 \begin{eqnarray}
 && f_{u^\bb,\mu}( x)-   f_{u^\bb,\mu}(0)
 \label{slid.11q}
 \\
 &&\qquad=\int_{-\ff}^\ff u^\bb(z-x) \,d\mu(z)-\int  _{-\ff}^\ff u^\bb( z) \,d\mu(z)
 \nonumber \\
 &&\qquad=\int _{-\ff}^\ff \(u^\bb(0)-u^\bb(z )\) -\(u^\bb(0)-u^\bb( z-x)\)  \,d\mu(z)
 \nonumber \\
 &&\qquad =\frac{1}{2}\int _{-\ff}^\ff \( (\si^\bb)^2(z )-(\si^\bb)^2(z-x)\) \,d\mu(z).
 \nonumber
 \end{eqnarray} 
 For $x \in [0,\De]$ it follows from (\ref{4.12w}) that,  
  \bea 
 2I_1(x) &:=&\int_{-\De}^\De \( (\si^\bb)^2(z )-(\si^\bb)^2(z-x)\) \,d\mu(z)\\
 &\le&(\si^\bb)^2(x) \int  _{-\De}^\De \frac{ 2u^\bb(0)}{u^\bb(z-x)+u^\bb(z)} \,d\mu(z)\nn .
\eea
 Consequently for $\De$   sufficiently small, 
\begin{equation} \label{4.21}
  \frac{ I_1(x ) }{(\si^\bb)^2(x)}\le C_{\De}\mu([-\De,\De]),
\end{equation}
for some $C_{\De}<\ff$ with $\lim_{\De\to 0}C_{\De}=2$.
  
  Now consider 
  \begin{equation} 
 2I_2(x ) :=\(\int _\De^\ff +\int _{-\ff}^{-\De}\)\( (\si^\bb)^2(z )-(\si^\bb)^2(z-x)\) \,d\mu(z).
\end{equation}
Note that for $x\in (0,\de)$ where $0<\de\le \De$,
  \begin{equation} 
  (\si^\bb)^2(z )-(\si^\bb)^2(z-x)=((\si^\bb)^2)'(z^*(z,x))x
\end{equation}
  for some $ z^*(z,x )\in[z-x,z ]$.   Therefore,
  \begin{equation} 
  \frac{ |I_2(x )|}{(\si^{\bb})^2(x)}= \frac{x}{(\si^{\bb})^2(x )}
\frac{1}{2} \(\int _\De^\ff +\int _{-\ff}^{-\De}\)((\si^\bb)^2)'(z^*(z,x))\,d\mu(z).\end{equation}
Since $|(\si^\bb)^2)'(z^*(z,,x))|$ is bounded on the range of integration we can use the Dominated Convergence Theorem to get,
 \begin{equation} \label{4.25}
 \lim_{x\downarrow 0} \frac{ |I_2(x )|}{(\si^{\bb})^2(x)}= \lim_{x\downarrow 0}\frac{x}{(\si^{\bb})^2(x )}
\frac{1}{2} \(\int _\De^\ff +\int _{-\ff}^{-\De}\)((\si^\bb)^2)'(z )\,d\mu(z).\end{equation}

When (\ref{1.20qq}) holds this is equal to 0. Using this  (\ref{newq}) and (\ref{4.21}) we see that in this case,  
 \begin{equation} 
  \lim_{x\downarrow 0} \frac{ f_{v^\bb,\mu}(x)}{(\si^{\bb})^2(x)}  \le C_{\De}\mu([-\De,\De]) + \frac{ f_{u^\bb,\mu}(0)}{2u^\bb(0,0)} .
\end{equation}
 Since this holds for all $\De>0$   and $\lim_{\De\to 0}C_{\De}=2$ we get (\ref{1.21qq}).

When  (\ref{1.22qq})  holds the right hand side of (\ref{4.25}) is equal to,
  \begin{equation} \label{4.27}
   \frac{1}{2C} \(\int _\De^\ff +\int _{-\ff}^{-\De}\)((\si^\bb)^2)'(z )\,d\mu(z),\end{equation}
which  in this case  gives,  
\begin{equation} \label{4.28}
  \lim_{x\downarrow 0} \frac{ f_{v^\bb,\mu}(x)}{(\si^{\bb})^2(x)}  \le   \frac{1}{2C}  \int _{-\ff}^\ff  ((\si^\bb)^2)'(z )\,d\mu(z)+ \frac{ f_{u^\bb,\mu}(0)}{2u^\bb(0,0)} .
\end{equation}
Using this we get (\ref{1.23qq}).\qed

 \begin{remark} {\rm  When  (\ref{1.22qq})  holds the integral in (\ref{4.28}) is   finite, since it, and the hypothesis that    $ (\si^\bb)^2 (t)\in C^1(T)$ and  $|(\si^\bb)^2 (t)'|$ is bounded on $(-\ff,-\De)\cup (\De,\ff)$ for all $ \De>0$ implies that $|((\si^\bb)^2 )'(t)|$ is bounded on $R^1.$  
 
 It is interesting to note that since $((\si^\bb)^2 )'(t)$ is an odd function whenever $\mu$ is a symmetric measure the integral in (\ref{4.28}) is equal to 0 and we get 
  \begin{equation} \label{1.2qmm}    \limsup_{  t \downarrow 0}\frac{   X_{\al } (t)}{  C  t \log\log 1/t }\le      1+ \frac{ f_{u^\bb,\mu}(0)}{2u^\bb(0)}    ,   \qquad a.s. 
   \end{equation}

}\end{remark}

 \begin{example} {\rm  In  \cite[1.25]{MRLIL} we point out that
 for $\bb>0$, 
  \begin{equation}      u^{\bb}_\ga(x-y) :=\frac{1}{ 2\pi}\int_{-\ff}^{\ff}\frac{ \cos\la (x-y)}{\bb+ \ga\la^2}\,d\la=\frac{1}{2 }\frac{e^{- ({ \bb/\ga})^{1/2}\,| x-y|} }{  (  \bb \ga )^{1/2}},\qquad x,y\in  R^{ 1}.  \label{exe.13ja}
\end{equation} 
In this example we simply take $\bb=K$ and $\ga=1/K$  so that,
   \begin{equation} \label{4.31}
  u_\ga^\bb(x-y)=\frac{1}{2}\exp^{-K|x-y|}.
\end{equation}
and
   \begin{equation} 
 ( \si_\ga^\bb)^2(x)=1- \exp^{-K|x|}.
\end{equation}

 Suppose  that $\mu$ is symmetric. Then
\begin{equation} 
  f_{u_\ga^\bb,\mu}(0)=  \frac{1}{2}\int_{-\ff}^\ff e^{-K|y|}\, \,d\mu(y)  = \int_{0}^\ff e^{-Ky}\, d\mu(y),
\end{equation}
 so that
\begin{equation}
\frac{ f_{u_\ga^\bb,\mu}(0)}{2u_ \ga^\bb(0)}= \int_{0} ^{\ff}e^{-Ky}\,d\mu(y).
\end{equation}
In addition,\begin{equation} \label{4.20qq}
 \lim_{x\da 0}\frac{(\si_\ga^\bb)  ^{2}(x)}{x}=  {K} .
\end{equation}
Therefore, it follows from (\ref{1.2qmm}) that,
\begin{equation} \label{4.26}
  \limsup_{  t \downarrow 0}\frac{   X_{\al,\bb } (t)}{K   t \log\log 1/t }\le   \(1+ \int_{0}^\ff e^{-Ky}\, d\mu(y)\),    \qquad a.s.  
  \end{equation}

Suppose  that $\mu$ is concentrated on the negative half line.  We evaluate the two terms on the right in (\ref{4.28}) for $u^\bb_\ga$ in (\ref{4.31}). Note that in this case $C$ in (\ref{4.28}) is equal to $K$. We have,
 \bea
  f_{u_\ga^\bb,\mu}  (0)&= &\int_{-\ff}^0 u_\ga^\bb(-y)\, d\mu(y)\nn\\
 &=& \frac{1}{2}\int_{-\ff}^0\exp^{-K|-y|}d\mu(y)=\frac{1}{2}\int_{-\ff}^0\exp^{Ky}d\mu(y).
\eea
Since $2u_ \ga^\bb(0)=1$ we get,
\be 
 \frac{ f_{u_\ga^\bb,\mu}  (0)}{2u_ \ga^\bb(0)}  =\frac{1}{2}\int_{-\ff}^0\exp^{Ky}d\mu(y).
\ee 
In addition,
\begin{eqnarray}
 && \frac{1}{2K}  \int _{-\ff}^0  ((\si^\bb)^2)'(z )\,d\mu(z)= \frac{1}{2K}  \int _{-\ff}^0  (1-\exp^{Kz})' \,d\mu(z)
 \label{}
 \\
 &&=-\frac{1}{2}\int_{-\ff}^0\exp^{Kz}d\mu(z).
 \nonumber
 \end{eqnarray}
Therefore, it follows from (\ref{1.2qmm}) and Theorem \ref{theo-7.2nn} that,
\begin{equation} \label{4.26q}
   \limsup_{  t \downarrow 0}\frac{   X_{\al,\bb } (t)}{K   t \log\log 1/t }=1,     \qquad a.s.  
   \end{equation}

}\end{example}

 \section{Proof of Theorem \ref{theo-1.5}}  
\ It is easy to check that, 
 \begin{equation} 
  u_{T_0}(x,x)+u_{T_0}(y,y)-2 u_{T_0}(x,y)=|s(x)-s(y)|,
\end{equation} 
and that
\begin{equation} 
   | u_{T_0}(x,v)- u_{T_0}(y,v)|\le |s(x)-s(y)|.
\end{equation}
Therefore,
\begin{equation} 
  |f_{u_{T_0},\mu}(x)-f_{u_{T_0},\mu}(y)|\le |s(x)-s(y)|\,\, |\mu|.
\end{equation}
The   condition that  $s\in C^2(0,\ff)$   shows that the criteria in Lemma \ref{lem-key}
are satisfied. The next lemma completes the proof of this theorem.

\begin{lemma} 
  \begin{equation} \label{5.4}
    f_{u_{T_0},\mu}(x)\sim s(x) |\mu|,\qquad \text{as $x\da 0$}.
\end{equation}
\end{lemma}

\Proof We have,
\begin{equation} 
   f_{u_{T_0},\mu}(x)=s(x)\int_{x}^\ff  \,d\mu(y)+ \int^{x} _0 s(y) \,d\mu(y).
\end{equation}
Therefore,
\begin{equation} \label{5.6}
  \frac{ f_{u_{T_0},\mu}(x)}{s(x)}=\int_{x}^\ff  \,d\mu(y)+ \int^{x} _0 \frac{s(y)}{s(x)} \,d\mu(y).
\end{equation} 
Since $s(y)$ is increasing   and $\mu$ is supported on $(0,\ff)$,  
\begin{equation}
\int^{x} _0 \frac{s(y)}{s(x)} \,d\mu(y)\leq \int^{x} _0  \,d\mu(y),
\end{equation}  Using this in (\ref{5.6}) we get (\ref{5.4}). 

The proof is completed as in the proofs of Theorems \ref{theo-1.4} and \ref{theo-1.5}.\qed

   \section{Proof of Theorem \ref{theo-1.6}} 
Let  
  \begin{equation} 
  (\wt \si^\bb)^2(x,y))=\wt u^\bb(x,x)+\wt u^\bb(y,y)-2\wt u^\bb(x,y) ,
\end{equation}
and define,    \begin{equation} \label{6.2qq}
  (\wt\si_v^\bb)^2(x,y))=\wt v^\bb(x,x)+\wt v^\bb(y,y)-2\wt v^\bb(x,y).
\end{equation}
  It follows from (\ref{diff.1v}) that
 \be \label{1.29sa}
 \wt v^{\bb} (x,y)= \wt u^\bb(x,y)-\frac{p_{\bb}(0)}{q_{\bb}(0)}q_{\bb}(x)q_{\bb}(y).    \ee 
Consequently,
\begin{equation} \label{3.xxqq}
  (\wt \si_v^\bb)^2(x,y)=(\wt \si^\bb)^2(x,y)-\frac{p_{\bb}(0)}{q_{\bb}(0)}\(q_{\bb}(x)-q_{\bb}(y) \)^2.
\end{equation}

  We show in  \cite[Lemma 4.9]{MRLIL} that 
   \be  \label{diff.14q}
  (\wt \si^\bb )^{2}  ( x, y)\sim
 \rho_{\bb}(0)|x-y|, \qquad \text{as $|x-y|\to 0$},
\ee
where, 
 \begin{equation} 
 \rho_{\bb}(0)= p_{\bb}'(0)q_{\bb}(0)-  q_{\bb}'(0)p_{\bb} (0) .
\end{equation} 
 
Since,  
\begin{equation} \label{6.6y}
  |q_{\bb}(x)- q_{\bb}(y)|\sim q_{\bb}'(0)|x-y|,\qquad \text{as $x,y\downarrow 0$},
\end{equation}
it follows   by (\ref{3.xxqq})   and (\ref{diff.14q})   that,  
\begin{equation} \label{3.xxww}
  (\wt \si_v^\bb)^2(x,y)\sim \rho_{\bb}(0) |x-y| ,\qquad \text{as $x,y\da 0$}. 
\end{equation}
By (\ref{diff.1v})
\begin{equation} \label{diff.1vqs}
\wt v^\bb (x,x)= \(p_{\bb}(x) -\displaystyle\frac{p_{\bb}(0)}{q_{\bb}(0)}q_{\bb}(x)\)q_{\bb}(x)  
 \ee  
Since,
\begin{equation} 
  p_{\bb}(x)=p_{\bb}(0)+p_{\bb}'(0)x+o(x)\quad \text{and}\quad
  q_{\bb}(x)=q_{\bb}(0)+q_{\bb}'(0)x+o(x), 
  \end{equation}  
as $x\da 0$, we see that, 
\bea \label{6.34s}
\wt v^\bb (x,x) &=&  \(\frac{p_{\bb}(x)}{p_{\bb}(0)} -\frac{q_{\bb}(x)}{q_{\bb}(0)}\)p_{\bb}(0)q_{\bb}(x) \nn\\
 &=&\(\frac{p_{\bb}'(0)}{p_{\bb}(0)} -\frac{q_{\bb}'(0)}{q_{\bb}(0)}\)p_{\bb}(0)q_{\bb}(x)x+o(x)\nn\\
&=&  {\rho_{\bb}(0)x} +o(x),\qquad\text{as $x\da 0$ } .
\eea

 Consequently,
 \begin{equation} \label{6.8}
 \wt v^\bb(x-y,x-y) \sim \rho_{\bb}(0) |x -y|,\quad \text{as $x,y\da 0$}, 
\end{equation}
and,  by (\ref{3.xxww})
\begin{equation} \label{6.10}
  (\wt \si_v^\bb)^2(x,y)\le C\wt  v^\bb(x-y,x-y) ,\quad \text{as $x,y\da 0$}, 
\end{equation}for some constant $C$.
This is the condition in (\ref{2.14}).

\medskip  Let  
  \begin{equation} \label{6.10q}
     f_{\wt u^\bb,\mu}(x)=\int_0^\ff \wt u^\bb(x,y) \,d\mu(y)
\end{equation}
and\begin{equation} \label{6.11q}
    f_{\wt v^\bb,\mu}(x)=\int_0^\ff  \wt v^\bb(x,y) \,d\mu(y).
\end{equation}
 Using (\ref{1.29sa}) we see that 
  \be \label{6.7}
  f_{\wt v^\bb,\mu}(x)
  = f_{\wt u^\bb,\mu}(x)-\frac{p_{\bb}(0)}{q_{\bb}(0)}q_{\bb}(x)\int_0^\ff  q_{\bb}(y)\,d\mu(y).
\ee
 Note that,\begin{equation} 
  f_{\wt u^\bb,\mu}(x)=p_{\bb}(x)\int_{x}^{\ff}q_{\bb}(y)\,d\mu(y)+q_{\bb}(x)\int_{0}^{x}p_{\bb}(y)\,d\mu(y)
\end{equation}
and,
\begin{equation} 
  f'_{\wt u^\bb,\mu}(x)=p_{\bb}'(x)\int_{x}^{\ff}q_{\bb}(y)\,d\mu(y)+q_{\bb}'(x)\int_{0}^{y}p_{\bb}(y)\,d\mu(y).
\end{equation}
 This shows that $ f_{\wt u^\bb,\mu}\in C^1([0,\ff))$   and  we see from (\ref{6.7}) that for $x,y\in[0,\de]$,
\be  
 \big| f_{\wt v^\bb, \mu}(x)-   f_{\wt v^\bb,\mu }(y)\big|
   \le   \sup_{s\in [0,\de]} f_{\wt u^\bb,\mu}'(s)|x-y|+\frac{p_{\bb}(0)}{q_{\bb}(0)}\big|q_{\bb}(x)-q_{\bb}(y)\big|\int_0^\ff  q_{\bb}(z)\,d\mu(z).  \ee
 Using (\ref{6.6y}) we see that,
 \be  
 \big|   f_{\wt v^\bb,\mu}(x)-    f_{\wt v^\bb, \mu}(y)\big|=O(|x-y|).
    \ee
 This is the  condition  in (\ref{2.15}).    
  
  \medskip We now consider the permanental process with kernel   \begin{equation} 
  v^\bb(t,t)=\wt v^\bb(t,t)+ f_{\wt v^\bb,\mu} (t).
\end{equation}
The behavior of $\wt v^\bb(t,t)$ at 0 is given in (\ref{6.8}). We now consider the behavior of  $f_{\wt v^\bb,\mu}$ at 0. 
 
%
%

  \begin{lemma}  
  \begin{equation} \label{6.23}
      f_{\wt v^\bb,\mu}(x) =  {x\,\rho_{\bb}(0)} \, \int_{0}^\ff \frac{q_{\bb}(y)}{q_{\bb}(0)} \,d\mu(y)  +o(x),\qquad \text{as $x \da 0$}.
\end{equation}

\end{lemma}

\Proof  
Rewrite (\ref{diff.1v}) as,  
 \begin{equation} \label{diff.1vq}
\wt v^\bb (x,y)= \left\{
 \begin{array} {cc}
 \(\displaystyle\frac{p_{\bb}(x)}{p_{\bb}(0)} - \frac{ q_{\bb}(x)}{q_{\bb}(0)}\)p_{\bb}(0) q_{\bb}(y),& \quad 0<x\leq y  
  \\
   \(\displaystyle\frac{p_{\bb}(y)}{p_{\bb}(0)} - \frac{ q_{\bb}(y)}{q_{\bb}(0)}\) p_{\bb}(0)q_{\bb}(x),& \quad 0<y\leq x.   
\end{array}  \right. 
\end{equation} 

  We have,  
\bea \label{6.20}
     f_{\wt v^\bb,\mu}(x)&=&\int^\ff_x \wt v^\bb(x,y)\,d\mu(y)+\int^x_0 \wt v^\bb(x,y) \,d\mu(y)\\
 &=&\(\frac{p_{\bb}(x)}{p_{\bb}(0)} -\frac{q_{\bb}(x)}{q_{\bb}(0)}\) p_{\bb}(0)\int_{x}^\ff q_{\bb}(y) \,d\mu(y)\nn\\
 &&\qquad+p_{\bb}(0)q_{\bb}(x)\int_{0}^x  \(\displaystyle\frac{p_{\bb}(y)}{p_{\bb}(0)} - \frac{ q_{\bb}(y)}{q_{\bb}(0)}\) \,d\mu(y)\nn\\
 &:=&I_{1}+I_{2}\nn. 
 \eea

It follows    as in (\ref{6.34s})   that,  
\bea \label{6.34}
   \frac{p_{\bb}(x)}{p_{\bb}(0)} -\frac{q_{\bb}(x)}{q_{\bb}(0)} &=& \(\frac{p_{\bb}'(0)}{p_{\bb}(0)} -\frac{q_{\bb}'(0)}{q_{\bb}(0)}\)x+o(x),\qquad\text{as $x\da 0$ }\\
&=& \frac{\rho_\bb(0)x}{p_{\bb}(0)q_{\bb}(0 )}+o(x),\qquad\text{as $x\da 0$ }\nn.
\eea
Using this and (\ref{6.20}) we see that,\begin{equation} 
  I_{1}(x)= x{\rho_{\bb}(0) }   \int_{0}^\ff \frac{q_{\bb}(y)}{q_{\bb}(0)} \,d\mu(y) +o(x),\qquad\text{as $x\da 0$. }
\end{equation}

We  complete the proof by showing that    $I_{2}(x)=o(x)$ as $x\da 0$.
%
Using (\ref{6.34}) we have,  
  \bea 
   I_{2}(x)&= & q_{\bb}(0)p_{\bb}(x)\int_{0}^x  \( \frac{p_{\bb}(y)}{p_{\bb}(0)}- \frac{ q_{\bb}(y)}{q_{\bb}(0)} \) \,d\mu(y)\\
   &= & q_{\bb}(0)p_{\bb}(x)\int_{0}^x  \(\frac{\rho_{\bb}(0)y}{p_{\bb}(0)q_{\bb}(0 )}+o(y) \) \,d\mu(y)\nn \\
     &= & \frac{ p_{\bb}(x)\rho_{\bb}(0) }{  p_{\bb}(0 )}\int_{0}^x y \,d\mu(y)+q_{\bb}(0)p(x)\int_{0}^x o(y)  \,d\mu(y)\nn .
\eea
The line above is,  
\begin{equation} 
  \le \frac{p(x)\rho_{\bb}(0)x }{ p_{\bb}(0 )}\int_{0}^x   \,d\mu(y)+q_{\bb}(0)p_{\bb}(x)o(x)\int_{0}^x   \,d\mu(y)=o(x ),\quad \text{as $x \da 0$},
\end{equation}
 since $\mu$ is supported on $(0,\ff)$.

  The proof is completed as in the proofs of Theorems \ref{theo-1.4} and \ref{theo-1.5}.\qed

  \begin{example}{\rm 
In   (\ref{diff.1}) we can take,
\begin{equation} 
  p_{\bb}(x)=e^{Kx},\qquad\text{and}\qquad  q_{\bb}(y)=\frac{1}{2}e^{-Ky};
\end{equation}
see  \cite[Example 9.3]{MRLIL}, in which case,
 \begin{equation} 
  \wt u^\bb (x,y)=\frac{1}{2}e^{-K|x-y|}.
\end{equation}It is easy to see that in this case,
\begin{equation} 
  \rho_{\bb}(0)= K.
\end{equation}
Consequently, (\ref{1.49}) gives,  
   \begin{equation} \label{opo}
   \limsup_{  t \downarrow 0}\frac{   X_{\al } (t)}{  Kt    \log\log 1/t }=  \(1+\int_{0}^\ff  e^{-Ky} d\mu(y)\),    \qquad a.s., 
   \end{equation}
   as in (\ref{4.26}).
}\end{example}

 \section{Additional observations}  

\begin{remark} {\rm \label{rem-7.1} In Theorems \ref{theo-1.3q} and \ref{theo-1.4} we require that   $(\si^0)^{2}$ and $(\si^\bb)^{2}$  are regularly varying at 0 with positive index. We show in  \cite[Section 8]{MRLIL} that these conditions are realized in many cases.  
}\end{remark}

\begin{remark}{\rm \label{rem-c1} We give a simple criteria for the conditions that
   $ (\si^\bb)^2 (t)\in C^1(0,\ff)$ and    $ |(\si^{\bb}) ^2 (x)'|$ is bounded on $[\de,\ff)$ for all $\de>0$.

\begin{lemma}   For all $\bb\geq 0$, if
\begin{equation} \label{7.3}
\sup_{0\le \la<\ff}   \frac{\la|\psi'(\la)|}{\bb+\psi(\la)} <\ff, 
\end{equation} 
   $  (\si^{\bb}) ^2 (x)\in C^1(0,\ff)$ and \be  \label{new4.58as}
 |(\si^{\bb}) ^2 (x)'|\le C  \frac{(\si^{\bb}) ^2 (x) }{x},\qquad \forall\,x\in(0,\ff),
\ee 
for  some constant $C $.

\end{lemma}

  \Proof    The proof is a modification of  \cite[Lemma 8.4]{MRLIL}   which gives sharper estimates for $ (\si^0)' (t) $ under  more restrictive conditions on $\psi(\la)$ .    
   When $x\neq 0$,  
\begin{equation} \label{4.59}
 (\si^{\bb}) ^2 (x) =\frac{4}{\pi }\int_{0}^{\ff}\frac{ \sin^2\la x/2}{ \bb+  \psi(\la)}\,d\la=\frac{4}{\pi x}\int_{0}^{\ff}\frac{ \sin^2t/2}{ \bb+   \psi(t/x)}\,dt,  
\end{equation}
 and,  
\bea \label{new4.58aw}
 (\si^{\bb}) ^2 (x)'&=&-\frac{  (\si^{\bb}) ^2 (x)}{x}+\frac{4}{\pi x}\frac{d}{dx}\(\int_{0}^{\ff}\sin^2 t/2 \frac{ 1}{ \bb+  \psi(t/x)}\,d t\)\nn\\
   &:=& -\frac{(\si^{\bb}) ^2 (x) }{x}+\nu_{\bb}(x). 
\eea
 We have,
 \begin{equation} 
  \nu_{\bb}(x)=\frac{4}{\pi x^2 } \int_{0}^{\ff} \frac{ \sin^2 t/2}{ \bb+   \psi(t/x)}\frac{(t/x)\psi'(t/x)}{  \bb+   \psi (t/x)  }\,dt.
  \end{equation}
  (See \cite[Lemma 8.3]{MRLIL} for the justification of  differentiating under the integral sign.)
Consequently, when  
\begin{equation} 
C:= \sup_{0\le \la<\ff}  \bigg|\frac{\la\psi'(\la)}{\bb+ \psi(\la)}\bigg|<\ff,
\end{equation}  
\begin{equation} \label{7.9}
    \nu_{\bb}(x)\le\frac{ C(\si^{\bb}) ^2 (x)}{  x  }.
\end{equation}
The inequalities in (\ref{new4.58as}) follow from (\ref{new4.58aw}) and (\ref{7.9}).\qed

Note that the condition in (\ref{7.3}) is satisfied for $\bb>0$ when $\psi(\la)$ is a normalized regularly varying function at  infinity, i.e.    for some $0<p\le 2$,
\begin{equation} 
  \psi(\la)=C\la^p\exp\(\int_{1}^\la\frac{\ep(u)}{u}\,du\),
\end{equation}
where $\lim_{u\to\ff}\ep(u)=0$.  When $\bb=0$ wecan also require that $\psi(\la)$ is a normalized regularly varying function at 0.

}\end{remark}

\begin{remark} \label{rem-7.2}{\rm 
  
 There are  kernels for  which the upper and lower bounds of the local moduli are the same.  

 \begin{theorem} \label{theo-7.2}   Let      $ Y $ be a strongly symmetric transient Borel right process with state space $T=R^{1}-\{0\}$ or $(0,\ff)$  and  continuous strictly positive  potential   densities $u(x,y)$ with respect to some $\si$--finite measure $m$ on $T$,   where    
    $ \lim_{t\downarrow0} u(t,t)=0.$
    
     Let $\{X_{\al,1 } (t),t\in [0,1] \}$ and $\{X_{\al,2 } (t),t\in [0,1] \}$ be an $\al-$permanental processes with kernels
 \be\label{7.00}
   v_{1}(s,t) =  u(s,t),\qquad \text{and }\qquad v_{2}(s,t)=u(s,t)+g(s)f(t),  \ee 
   respectively,  where  $f$ is given in (\ref{rp.1})  and  
     \begin{equation} \label{1.17mmaa}
  g(s)=\int_{T} u(s,v)\,d\nu(v)\end{equation}
         for some positive finite measure $\nu$. 
If $u(v,v)$
is   regularly varying   at zero   with positive index and 
\begin{equation} \label{2.14m}
   u ( s,s) +    u( t,t)-    2u ( s,t) \le C u(|t-s|,|t-s|).  
\end{equation} 
then
\begin{equation} \label{6.2ssq1}
  \limsup_{t\da 0}\frac{X_{\al,1}(t )}{  u(t,t)\log \log 1/t }= 1 . \end{equation}
\et
If, in addition, 
 \begin{equation} \label{2.15m}
  \big|f(s)-    f(t)\big|\le C u(|t-s|,|t-s|)\quad \text{and}\quad \big|g(s)-    g(t)\big|\le C u(|t-s|,|t-s|) 
\end{equation} 
for some constant $C$ for all $s,t\in [0,\de]$, for some $\de>0$,
then    (\ref{6.2ssq1}) holds for $X_{\al,2}$.
 (The existence of     $X_{\al,2}$ follows from \cite[Theorem 1.1]{MRLIL}.) \qed 
 
 \medskip \noindent {\bf Proof of Theorem \ref{theo-7.2} } When the kernel is $v_1(s,t)$ the lower bound follows from Theorem \ref{theo-7.2nn} and the upper bound follows from Theorem \ref{theo-2.2}.

Suppose the kernel is $v_2(s,t)$.
Consider  $ \wt  \si   ^2(s,t)$. Using   (\ref{7.00}) and the fact that $u(s,t)$ is symmetric we can
  write this as, 
 \bea \label{3.107qqsm}
 \wt   \si  ^2(s,t)  &=&   v_2 ( s,s) +    v_2( t,t)-    v_2 ( s,t) -   v_2( t,s)\\
 &&\qquad    + \(v_2 ^{1/2}( s,t)-       v_2^{1/2} ( t,s)\)^2\nn\\&=& u ( s,s) +    u( t,t)-    2u ( s,t) +(  f(s)- f(t))( g(s)-  g(t))\nn\\
  &&\qquad    + \(v_2^{1/2}(s,t)- v_2^{1/2}(t,s)\)^2.\nn \eea
 Using   (\ref{2.14m}) and (\ref{2.15m})  we see  that,
 \begin{eqnarray}
 && u ( s,s) +    u( t,t)-    2u ( s,t) +(  f(s)- f(t))( g(s)-  g(t))
 \label{3.107qq1}
 \\ 
 &&\qquad \leq C'u(|t-s|,|t-s|),
 \nonumber
 \end{eqnarray}
  for some constant $C'$, for all $s,t\in [0,\de]$. Using the fact that $u(s,t)$ is symmetric we   have,
 \bea \label{3.107qq2w}
 && \(v_2^{1/2}(s,t)- v_2^{1/2}(t,s)\)^2 \le   \big| v_2 (s,t)- v_2 (t,s)\big|\\ &&  \qquad=\big|f(s)g(t)-   f(t)g(s)\big|=\big|f(s)(g(t)-g(s))-   (f(t)-f(s))g(s)\big|\nn\\
 && \qquad \le f(s)\big|g(s)-    g(t)\big|+g(s)\big|f(s)-    f(t)\big|\nn.
\eea
  Consequently,  
\be  \label{3.107qq}
    \si ^2(s,t)  \le   C' u(|t-s|,|t-s|) +2\( f(s)\big|g(s)-    g(t)\big|+g(t)\big|f(s)-    f(t)\big|\),\ee 
and by (\ref{1.17mms}), 
\bea 
 &&   f(s)\big|g(s)-    g(t)\big|+g(t)\big|f(s)-    f(t)\big|\\
 &&\qquad \le  |\mu|u(s,s)\big|g(s)-    g(t)\big|+  |\nu|u(t,t)\big|f(s)-    f(t)\big|.\nn
\eea
It now follows by a minor generalization of Theorem \ref{theo-2.2} that,
\begin{equation} \label{6.2ss}
  \limsup_{t\da 0}\frac{X_{\al,2}(t )}{( u(t,t)+C'u^2(t,t))\log \log 1/t }\le 1 ,\qquad \text{a.s.}, \end{equation}
for some constant $C'$. This implies that
\begin{equation} \label{6.2ssq}
  \limsup_{t\da 0}\frac{X_{\al,2}(t )}{  u(t,t) \log \log 1/t }\le 1 ,\qquad \text{a.s.} \end{equation}
  
The lower bound for $X_{\al,2}$ follows from a minor modification of the proof of Theorem \ref{theo-7.2nn}.   
 Let 
  $U_{f,g}$ be a non-singular $n\times n$ matrix,   \bea \hspace{-.3 in}  \left (
\begin{array}{ cccc } 1 &f(1 )&\ldots&f(n )  \\
g(1)  &U_{1 ,1 }+g(1)f(1 )&\ldots&U_{1 ,n }+g(1)f(n )  \\
\vdots& \vdots &\ddots &\vdots  \\
g(n)   &U_{n ,1 }+g(n)f(1 )&\ldots&U_{n ,n }+g(n)f(n )  
\end{array}\right ).\label{19.39q}
   \eea 
where $g$ is any function. 

  Let  $ U_{f,g}^{-1}$   denote the inverse of $U_{f,g}$ and $U_{f,g}^{j,k}$   denote the elements of $U_{f,g}^{-1}$.  
One can check that, 
\be  U_{f,g}^{-1} =\left (
\begin{array}{ cccc  }  1+\rho' &-\sum_{i=1}^{n}f(i)U ^{i ,1 }   &\dots&-\sum_{i=1}^{n}f(i)U ^{i ,n }  \\
- \sum_{j=1}^{n}g(1)U ^{1 ,j } & U ^{1 ,1 } & \dots &  U ^{ 1 ,n }    \\
\vdots&\vdots& \ddots&\vdots  \\
- \sum_{j=1}^{n}g(n)U ^{n ,j }& U ^{n ,1 }&  \dots & U ^{n ,n }   \end{array}\right ),\label{19.40}
   \ee 
   where 
   \begin{equation}
\rho' =\sum_{j=1}^{n}\sum_{i=1}^{n}g(j)f(i)U^{i ,j} .\label{19.10v}
\end{equation}
As in (\ref{7.7}) we have 
 \begin{equation} \label{7.7qq}
  U _{f,g} ^{i,i}=  U^{i-1,i-1},\qquad i=2,\ldots,n+1.
\end{equation}
The proof now proceeds like the proof of Theorem  \ref{theo-7.2nn}.\qed

}\end{remark}


\begin{thebibliography} {10}     

	 

  
\bibitem {EK}
N. Eisenbaum  and  H. Kaspi,
\newblock On permanental processes, {\em Stochastic Processes and their Applications},
{\em 119},  (2009),  1401-1415.

  \bibitem{Breiman} L.  Breiman, {\em Probability}, Classics  in Applied Mathematics, SIAM, Philadelphia,  (1992).
 

\bibitem{book} M. B.  Marcus and J.~Rosen, {\em Markov Processes,
Gaussian Processes and Local Times}, Cambridge University Press, New
York,  (2006).
 
%
%
%
%

\bibitem{MRejp} M. B. Marcus and J.~Rosen,
Sample path properties of permanental processes,  {\it Electronic Journal of Probability},\, Volume 23 (2018), no. 58, 1-47.

\bibitem{MRnonsym} M. B. Marcus and J.~Rosen, Permanental processes with kernels that are not equivalent to a symmetric matrix.   High Dimensional Probability VIII, The Oaxaca Volume, Progress in Probability 74,\,(2019), 305-320. Birkhauser, Boston.

\bibitem{MRLIL} M. B. Marcus and J.~Rosen,
 Law of the iterated logarithm for $k/2$--permanental processes and the local times of   related Markov processes, Submitted. 
 
 \bibitem{MRun} M. B. Marcus and J.~Rosen,
 Conditions for permanental processes to be unbounded, ... 

    
%



%
%


 \bibitem{RY}
Daniel Revuz and Marc Yor.
\newblock {\em Continuous martingales and {B}rownian motion}, volume 293.
\newblock Springer-Verlag, Berlin, third edition, 1999.
 
 \end{thebibliography}
 \end{document}